\documentclass[12pt]{article}
\usepackage{amsfonts,amsmath,amscd}
\usepackage[a4paper,left=3cm,right=2cm,top=3cm,bottom=2cm]{geometry}
\usepackage{hyperref}
\usepackage{lmodern}
\usepackage[utf8x]{inputenx} %
\usepackage[LGR, T1]{fontenc}
\usepackage{textcomp} 
\usepackage{graphicx}
\usepackage{wrapfig}

\newcommand{\R}{\mathbb{R}}

\newtheorem{defi}{Definition}[section]
\newtheorem{obs}{Remark}[section]

\newtheorem{exam}{Example}[section]
\newtheorem{teorema}{Theorem}[section]
\newtheorem{coro}{Corollary}[section]

\newcommand{\n}{\noindent}

\begin{document}

\title{Invariant Solutions for Gradient Ricci Almost Solitons}

\author{
	\textbf{Benedito Leandro}
	\\
	{\small\it CIEXA-Universidade Federal de Jata\'i }\\
	{\small\it BR 364, km 195, 3800, 75801-615,
		Jata\'i, GO, Brazil. }
	\\
	{\small\it e-mail:  bleandroneto@gmail.com}
	\\
	\textbf{Romildo Pina \footnote{ Partially supported by
			CAPES-PROCAD.}}
	\\
	{\small\it IME, Universidade Federal de Goi\'as,}\\
	{\small\it  131, 74001-970, Goi\^ania, GO, Brazil. }\\
	{\small\it e-mail: romildo@ufg.br }
	\\
\textbf{Tatiana Pires Fleury Bezerra}
\\
{\small\it IFG -Instituto Federal de Educa\c{c}\~{a}o, Ci\^{e}ncia e Tecnologia de Goi\'as }\\
{\small\it 74968-755, Lt-1A, Parque Itatiaia, Ap. de Goi\^ania, GO,
Brazil.}
\\
{\small\it e-mail:  tatipibe@hotmail.com}
}

\maketitle
\thispagestyle{empty}

\markboth{abstract}{abstract}
\addcontentsline{toc}{chapter}{abstract}

\begin{abstract}

     \noindent
 In this paper we provide an ansatz that reduces a pseudo-Riemannian gradient Ricci almost soliton (PDE) into an integrable system of ODE. First, considering a warped structure with conformally flat base invariant under the action of an $(n-1)$-dimensional translation group and semi-Riemannian Einstein fiber, we provide the ODE system which characterizes all such solitons. Then, we also provide a classification for a conformally flat pseudo-Riemannian gradient Ricci almost soliton invariant by the actions of a translation group or a pseudo-orthogonal group. Finally, we conclude with some explicit examples.
         
\end{abstract}


\section{Introduction and main statements}

In the early eighties, Jean-Pierre Bourguignon introduced a flow to study the evolution of the Ricci curvature and the Einstein metrics. The {\it Ricci Bourguignon flow} is given by
$$\frac{\partial}{\partial\,t}g(t)=-2(Ric-\kappa\,R\,g)(t),$$ where $Ric$ and $R$ are, respectively, the Ricci tensor and the scalar curvature for the metric $g$. This flow is an interpolation between the Ricci flow and the Yamabe flow (see \cite{catino} and the references therein). 

The self-similar solutions of this flow are called {\it Einstein solitons} and correspond to the
equation $$Ric+Hess(h)=(\kappa\,R+\mu\,)g,$$ where $\kappa$, $\mu\in\mathbb{R}$. 
When we replace in the above equation the term $\kappa R+\mu$  for $\lambda$, an arbitrary smooth function, we can call this equation a {\it gradient Ricci almost soliton}. We say that a {\it Ricci almost soliton} is a smooth manifold satisfying $$Ric + \mathcal{L}_{X}g=\lambda g,$$ where $\mathcal{L}_{X}g$ represents the Lie derivative of the metric $g$ with respect to a tangent vector field $X$ and $\lambda$ is an arbitrary smooth function. The Ricci almost solitons generalize Einstein solitons. They are also natural generalizations for Ricci solitons, which are self-similar solutions for the Ricci flow.

\begin{defi}
A smooth manifold $(M^{n}, g)$ is a gradient Ricci almost soliton if there exist two smooth functions $h$ and $\lambda$ on $M$ such that
\begin{eqnarray}\label{eq1}
Ric_{g}+Hess_{g}(h)=\lambda g,
\end{eqnarray}
where Ric$_g$ is the Ricci tensor, Hess$_g(h)$ is the Hessian of the potential function $h$ with respect to the metric $g$.
\end{defi}

A smooth manifold $(M^{n}, g)$ is a \emph{ \it gradient Ricci soliton} if there exist a smooth function $h:M\longrightarrow\mathbb{R}$ (called the potential function) and a constant $\lambda$ satisfying \ref{eq1}. A gradient Ricci soliton is said to be shrinking, steady, or expanding if $\lambda>0$, $\lambda=0$, or $\lambda<0$, respectively. When $M$ is a Riemannian manifold, usually, one requires the manifold to be complete. In the case of semi-Riemannian manifolds, one does not require $(M,g)$ to be complete (see \cite{BBGG},   \cite{BCGG}, \cite{BGG} e  \cite{Onda}).

We observe that \ref{eq1} can be considered a perturbation of the Einstein equation
$$\mbox{Ric}_g=\rho g, \qquad \rho\in\R.$$
When $h$ is constant, we call the underlying Einstein manifold a trivial Ricci soliton.

Pigola et al. \cite{PRRS}, explored the Ricci almost solitons in a very comprehensive way. They provided topological properties, volume comparison results, a gap theorem and some explicit examples are given. 

Since gradient Ricci almost solitons contain gradient Ricci solitons as a particular case, we can say that the gradient Ricci almost soliton is \emph{ \it proper} if the function $\lambda$ is non-constant (see \cite{BGRR}). 

In \cite{FFGP} the authors presented a necessary and sufficient condition for constructing gradient Ricci almost solitons that are realized as warped products. They provided an example of a particular Riemannian solution of the PDEs that arise from the hypothesis that the base is conformally flat and invariant by translation, in which the fiber is an Einstein manifold. Here, we classify all pseudo-Riemannian gradient Ricci almost soliton with warped structure in Theorem \ref{teo1.3}. Moreover, Theorem \ref{coro1.3} provides a method capable of producing an infinite number of pseudo-Euclidean warped product gradient Ricci almost soliton such that the base is invariant under the action of an $(n-1)$-translation group and the fiber is Ricci flat (see explicit examples below).
 
Further, Barros, Batista and Ribeiro Jr \cite{BBR}, proved that either a Euclidean space $\mathbb{R}^{n}$ or a standard sphere $\mathbb{S}^{n}$ is the unique manifold with non-negative scalar curvature, which carries a structure of a gradient Ricci almost soliton, provided this gradient is a nontrivial conformal vector field. Also, they showed that a compact locally conformally flat almost Ricci soliton is isometric to a Euclidean sphere $\mathbb{S}^{n}$ provided that an integral condition holds. In addition, they constructed examples of gradient Ricci almost solitons. Consider the warped product manifold $M^{n+1}=\mathbb{R}\times_{\cosh(t)}\mathbb{S}^{n}$ with metric $g=dt^{2}+\cosh^{2}(t)g_{0}$, where $g_{0}$ is the standard metric of $\mathbb{S}^{n}$. They proved that $(M^{n+1}, g, \nabla h, \rho )$, where $h(x, t)=\sinh(t)$ and $\lambda(x,t)=\sinh(t) +n$, is a gradient Ricci almost soliton.

 Barros et al. \cite{BGR} proved that a gradient Ricci almost soliton $(M^{n}, g,\nabla h, \lambda )$, whose Ricci tensor is Codazzi, has constant sectional curvature. In particular, in the compact case, they deduced that $(M^{n}, g)$ is isometric to a Euclidean sphere and $h$ is a height function. They classified gradient Ricci almost solitons with constant scalar curvature and provided a suitable function that achieved a maximum in $M^{n}$.

In 2011, Catino \cite{CA} introduced the notion of generalized quasiâ€“Einstein manifold that generalizes the concepts of Ricci soliton, Ricci almost soliton, and quasiâ€“Einstein manifolds. He proved that a complete generalized quasiâ€“Einstein manifold with harmonic Weyl tensor and zero radial Weyl curvature is locally a warped product with $(n-1)$-dimensional Einstein fibers. Furthermore, in this paper, Catino proved the following result: 

    ``Let $(M^{n}, g)$, $n\geq3$ be a locally conformally flat gradient Ricci almost soliton. Then, around any regular point of $h$, the manifold $(M^{n}, g)$ is locally a warped product with $(n-1)$-dimensional Einstein fibers of constant sectional curvature.''
    
 In particular, this implies a local characterization for locally conformally flat gradient Ricci almost solitons, similar to that proved for gradient Ricci solitons. 

 The local structure of half conformally flat gradient Ricci almost solitons was recently investigated in \cite{BGRR}, showing they are locally conformally flat in a neighborhood of any point where the gradient of the potential function is non-null. In \cite{CEGGV}, the authors proved that a locally homogeneous proper Ricci almost soliton has either
constant sectional curvature or is locally isometric to a product $ \mathbb{R} \times N(c)$, where $N(c)$ is a space of constant curvature.
 
Inspired by the local classification of conformally flat gradient Ricci almost soliton that we discussed above, we consider a warped product structure for gradient Ricci almost solitons which is not necessarily conformally flat and then we classify such solitons.

Considering $(B, g_{B})$ and $(F, g_{F})$ semi-Riemannian manifolds, with $f>0$ being a smooth function on the base $B$, the warped product $M=B\times_{f} F$, with {\it fiber} $F$ and {\it warping function} $f$, is the product manifold $M=B\times F$ furnished with metric tensor
              
$$ \tilde{g}=g_{B}+f^{2}g_{F}.$$

In what follows, we find a family of gradient Ricci almost solitons in the case of the warped product $(M, \tilde{g})=(\mathbb{R}^{n},\bar{g})\times_{f}(F^{m}, g_{F})$, where the fiber is a semi-Riemannian Einstein manifold and the base is conformal to a pseudo-Euclidean space which is invariant under the action of an $(n-1)$-dimensional translation group.

 We denote $\psi,_{x_{i}}$, $f,_{x_{i}}$, and $h,_{x_{i}}$, first order derivative and $\psi,_{x_{i}x_{j}}$, $f,_{x_{i}x_{j}}$, and $h,_{x_{i}x_{j}}$ as the second order derivative of functions $\psi$, $f$, and $h$, with respect to $x_{i}$ and $x_{i}x_{j}$, respectively. Moreover, for any smooth function $W$ we denote $$W^{'}=\frac{dW}{d\xi}\quad\mbox{or}\quad W^{'}=\frac{dW}{dr}.$$
Without further ado, we state our main results.

\begin{teorema} \label{teo1.3} 	    Let $(\mathbb{R}^{n}, g)$ be a pseudo-Euclidean space with Cartesian coordinates $x=(x_{1}, ..., x_{n})$, $g_{ij}=\delta_{ij}\varepsilon_{i}$, where $\delta_{ij}$ is the delta
Kronecker, $\varepsilon_{i}=\pm1$ with at least one $\varepsilon_{i}=1$. Consider $(M, \tilde{g})=(\mathbb{R}^{n},\bar{g})\times_{f}(F^{m}, g_{F}),$ where $\overline{g}=\frac{1}{\psi^{2}}g$ and $F$ is a semi-Riemannian Einstein manifold with constant Ricci curvature $\lambda_{F}$. Moreover, assume non-constant smooth functions $h(\xi)$, $\lambda(\xi)$ and $f(\xi)>0$, where $\xi=\sum^{n}_{i=1}\alpha_{i} x_{i},\; \alpha_{i} \ \in \mathbb{R} $ and $\sum^{n}_{i=1} \varepsilon_{i}\alpha_{i}^{2}=\varepsilon_{i_{0}}$ equals $-1$ or $1$ if  is a timelike or spacelike
vector, respectively. Then, the warped product metric $\tilde{g}=\bar{g}+f^{2}g_{F}$ is a gradient Ricci almost soliton \ref{eq1} with $h$ as a potential function if, and only if, the functions $\psi$, $f$, $\lambda$, and $h$ satisfy
 	    \begin{equation} \label{eq7}
 	             \left\{ \begin{array}{lcl}
 	                 f\left[ (n-2)\psi^{''}+2\psi^{'}h^{'}+\psi h^{''}\right]-m\psi f^{''}-2m \psi^{'}f^{'} =0;\\ \\
 	                 
 	                 \varepsilon_{i_{0}} \left[f \psi\psi^{''} -(n-1)f\left( \psi^{'} \right)^{2} +m\psi\psi^{'}f^{'}-f\psi\psi^{'}h^{'} \right]  = \lambda f;\\ \\
 	                 
 	                 \varepsilon_{i_{0}} \left[-f \psi^{2} f^{''}+(n-2)f \psi f^{'}\psi^{'}-(m-1) \psi^{2} \left( f^{'}\right)^{2}+f \psi^{2}f^{'}h^{'} \right]  = \lambda f^{2}-\lambda_{F}.   
 	              \end{array} \right.
 	       \end{equation}  
\end{teorema}

In the next result we prove that if $f\psi=1$ and $F$ is Ricci flat, then the metrics $\tilde{g}$ are gradient Ricci almost solitons.

\begin{teorema} \label{coro1.3}
    	Let $(\mathbb{R}^{n}, g)$ be a pseudo-Euclidean space with coordinates $x=(x_{1}, ..., x_{n})$, $g_{ij}=\delta_{ij}\varepsilon_{i}$, where $\delta_{ij}$ is the delta
Kronecker, $\varepsilon_{i}=\pm1$ with at least one $\varepsilon_{i}=1$. Consider $M=(\mathbb{R}^{n},\bar{g})\times_{f} F^{m},$ a warped product, where $\overline{g}=\frac{1}{\psi^{2^{}}}g,$ $F$ is a Ricci-flat semi-Riemannian Einstein manifold. Consider $f(\xi)>0$, $\lambda(\xi)$ e $h(\xi)$ non-constant smooth functions, where $\xi=\sum^{n}_{i=1}\alpha_{i} x_{i},\ \ \alpha_{i} \ \in \mathbb{R}$ and $\sum^{n}_{i=1} \varepsilon_{i}\alpha_{i}^{2}=\varepsilon_{i_{0}}$ equals $-1$ or $1$ if  is a timelike or spacelike
vector, respectively. Given any function $ \psi (\xi)$, the warped product metric $\tilde{g}=\bar{g}+f^{2}g_{F}$ is a gradient Ricci almost soliton \ref{eq1} with $h$ as a potential function, where the functions $f$, $h$, and $\lambda$ are given by
    	
    	       \begin{equation} \label{eq8}
    	                \left\{ \begin{array}{lcl}
    	                      f(\xi)\psi(\xi)=1;\\ \\
    	       
    	                      h(\xi) = k+ \displaystyle \int \left\lbrace c-(m+n-2) \displaystyle \int \psi \psi^{''}  d\xi \right\rbrace \frac{1}{\psi^{2}} d\xi;\\ \\
    	       
    	                      \lambda(\xi)=\varepsilon_{i_{0}} \left\lbrace \psi\psi^{''}-(m+n-1)(\psi^{'})^{2}-c \dfrac{\psi^{'}}{\psi}+ (m+n-2) \dfrac{\psi^{'}}{\psi} \displaystyle\int \psi\psi^{''}d\xi \right\rbrace, 
    	                   \end{array} \right.
    	       \end{equation}		
    	 
     where $c$ and $k$ are constants.
\end{teorema}
    
\begin{obs}\label{mas}
In Theorem \ref{coro1.3} $f\psi=1$, thus the metric $\tilde{g}$ can be rewritten as

     $$\tilde{g}=\bar{g}+f^{2}g_{F}=\frac{1}{\psi^{2^{}}}g_{E}+\left(\frac{1}{\psi} \right) ^{2}g_{F}=\frac{1}{\psi^{2^{}}}  \left(g_{E}+g_{F} \right).$$    
               
 Thus, all metrics conformal to the product manifold $\left( R^{n}\times F^{m}\right)$, invariant by translation, where $F$ is a Ricci flat manifold, are gradient Ricci almost solitons.
 
\end{obs}

Hereafter, we find a family of gradient Ricci almost solitons in the case of $(M, \tilde{g})=(\mathbb{R}^{n},\bar{g})$, where $M$ is conformal to a pseudo-Euclidean space which is invariant under the action of a pseudo-orthogonal group.

Let $(\mathbb{R}^{n}, g)$ be the standard pseudo-Euclidean space with metric $g$ and coordinates $x=(x_{1}, ..., x_{n})$ with $g_{ij}=\delta_{ij}\varepsilon_{i}, 1\leq i, j \leq n$, where $\delta_{ij}$ is the delta
Kronecker, $\varepsilon_{i}=\pm1$ with at least one $\varepsilon_{i}=1$. Let $r=\displaystyle\sum^{n}_{i=1}\varepsilon_{i} x^{2}_{i},\ \ \alpha_{i} \ \in \mathbb{R},$ be a basic invariant for an $(n-1)-$dimensional pseudo-orthogonal group. We want to obtain differentiable functions $f(r)$, $\psi(r)$ and $\lambda(r)$ such that the metric $\overline{g}=\frac{1}{\psi^{2}}g $ satisfies \ref{eq1}.

We  show that in the Riemannian case, all metrics conformal to a Euclidean metric and invariant by rotation are gradient Ricci almost solitons, see Corollary \ref{coro1.1}. Moreover, in Corollary \ref{coro1.1}, given any function $\psi(r)$, there are $f(r)$ and $\lambda(r)$ such that the metric $\overline{g}$ is a gradient Ricci almost soliton. This provides a method to build many examples of solitons invariants by rotation.

In Corollary \ref{coro1.1} and Theorem \ref{teo1.1}, we consider metrics conformal to pseudo-Euclidean spaces, then we find families of gradient Ricci almost soliton invariants by a pseudo-orthogonal group action.

\begin{teorema} \label{teo1.1}
	Let $(\mathbb{R}^{n}, g)$ be a Euclidean space $n\geq3$ with coordinates $x=(x_{1}, ..., x_{n})$, $g_{ij}=\delta_{ij}\varepsilon_{i}$, where $\delta_{ij}$ is the delta
Kronecker, $\varepsilon_{i}=\pm1$ with at least one $\varepsilon_{i}=1$. Consider non-constant smooth functions $h(r)$ and $\psi(r)$, where $r=\sum^{n}_{i=1}\varepsilon_{i}x_{i}^{2}$. There exists metric $\overline{g}=\frac{1}{\psi^{2}}g $ such that $(\mathbb{R}^{n}, \overline{g})$ is a gradient Ricci almost soliton \ref{eq1} with $h$ as a potential function if, and only if, the functions $h$, $\psi$ and $\lambda$ satisfy    	    
    	\begin{equation} \label{eq3}
	        \left\{ \begin{array}{lc}
	            (n-2)\psi^{''}+2\psi^{'}h^{'}+\psi h^{''} =0;\\ \\
	
	            4(n-1) \psi\psi^{'}+4r \psi\psi^{''} -4(n-1)r\left( \psi^{'}\right) ^{2} -4r\psi\psi^{'} h^{'} +2\psi^{2}h^{'}=\lambda.   
	        \end{array} \right.
	\end{equation}
\end{teorema}

Note that, the conformal function is free for choice. Therefore, if we choose a conformal function for Theorem \ref{teo1.1}, we can build a pseudo-Riemannian gradient Ricci almost soliton invariant by the action of a pseudo-orthogonal group (rotational in the Riemannian case), provided that the system \ref{eq3} is integrable. We sum up this discussion in the following corollary.

\begin{coro} \label{coro1.1}
	Let $(\mathbb{R}^{n}, g)$ be a Euclidean space $n\geq3$ with coordinates $x=(x_{1}, ..., x_{n})$, $g_{ij}=\delta_{ij}\varepsilon_{i}$, where $\delta_{ij}$ is the delta
Kronecker, $\varepsilon_{i}=\pm1$ with at least one $\varepsilon_{i}=1$. Consider non-constant smooth functions $h(r)$ and $\psi(r)$, where $r=\sum^{n}_{i=1}\varepsilon_{i}x_{i}^{2}$. Given any function $ \psi(r)$, the metric $ \overline{g}=\frac{1}{\psi^{2^{}}}g $ is a gradient Ricci almost soliton \ref{eq1} with $h$ as a potential function, where the functions $h$ and $\lambda$ are given by 
	
	      \begin{equation} \label{eq4}
            	\left\{ \begin{array}{llc}
	                h(r)&=&\displaystyle\int \left[ c-(n-2)\displaystyle\int\psi\psi^{''}dr \right] \dfrac{1}{\psi^{2}}dr +k;\\ \\
	
\lambda(r)&=& 4(n-1) \psi\psi^{'}+4r \psi\psi^{''} -4(n-1)r\left( \psi^{'}\right) ^{2}\\ \\
 &-&4cr\frac{\psi^{'}}{\psi} -2(n-2) \left( 1-2r \frac{\psi^{'}}{\psi} \right) \displaystyle\int\psi\psi^{''}dr +2c,   
             	\end{array} \right.
	      \end{equation}
	where $c$ and $k$ are constants.
\end{coro}

For our next results, let $(\mathbb{R}^{n}, g)$ be the standard pseudo-Euclidean space with metric $g$ and coordinates $x=(x_{1}, ..., x_{n})$ with $g_{ij}=\delta_{ij}\varepsilon_{i}, 1\leq i, j \leq n$, where $\delta_{ij}$ is the delta
Kronecker, $\varepsilon_{i}=\pm1$ with at least one $\varepsilon_{i}=1$. Let $\xi=\displaystyle\sum^{n}_{i=1}\alpha_{i} x_{i},\ \ \alpha_{i} \ \in \mathbb{R},$ be a basic invariant for an $(n-1)-$dimensional pseudo-orthogonal group, where $\sum^{n}_{i=1} \varepsilon_{i}\alpha_{i}^{2}=\varepsilon_{i_{0}}$ equals $-1$, $0$, or $1$ if it is a timelike, lightlike, or spacelike
vector, respectively. We want to obtain differentiable functions $f(r)$, $\psi(r)$ and $\lambda(r)$ such that the metric $\overline{g}=\frac{1}{\psi^{2}}g $ satisfies \ref{eq1}.

Now, we first obtain the necessary and sufficient conditions on $f(\xi)$ and $\psi(\xi)$ for the existence of $ \overline{g}$. These conditions differ depending on the direction $\alpha=\displaystyle\sum^{n}_{i=1}\alpha_{i} \dfrac{\partial}{\partial x_{i}}$ being timelike or spacelike. Remember that we are considering proper solitons.

\begin{teorema} \label{teo1.2}
    		Let $(\mathbb{R}^{n}, g)$ be a pseudo-Euclidean space $n\geq3$ with coordinates $x=(x_{1}, ..., x_{n})$, $g_{ij}=\delta_{ij}\varepsilon_{i}$. Consider non-constant smooth functions $h(\xi)$ and $\psi(\xi)$, where $ \xi=\displaystyle\sum^{n}_{i=1}\alpha_{i} x_{i},\ \ \alpha_{i} \ \in \mathbb{R} $ and $\displaystyle \sum^{n}_{i=1} \varepsilon_{i}\alpha_{i}^{2}=\varepsilon_{i_{0}}\neq 0$. There exists the metric $ \overline{g}=\frac{1}{\psi^{2^{}}}g $ such that $(\mathbb{R}^{n}, \overline{g})$ is a gradient Ricci almost soliton \ref{eq1} with $h$ as a potential function if, and only if, the functions $h$, $\psi$, and $\lambda$ satisfy
    	
    	 \begin{equation} \label{eq5}
    	        \left\{ \begin{array}{lc}
    	                 (n-2)\psi^{''}+2\psi^{'}h^{'}+\psi h^{''} =0;\\ \\
    	
                         \varepsilon_{i_{0}} \left[ \psi\psi^{''}-(n-1) \left( \psi^{'} \right)^{2} - \psi\psi^{'} h^{'} \right]= \lambda.
             	\end{array} \right.
    	\end{equation}
    	
\end{teorema}
 
In the next result we provide families of gradient Ricci almost solitons which are invariant under the action of an $(n-1)$-dimensional translation group.  

\begin{coro} \label{coro1.2}
    		Let $(\mathbb{R}^{n}, g)$ be a pseudo-Euclidean space $n\geq3$ with coordinates $x=(x_{1}, ..., x_{n})$, $g_{ij}=\delta_{ij}\varepsilon_{i}$. Given any function $ \psi(\xi)$, the metric $ \overline{g}=\frac{1}{\psi^{2^{}}}g $ is a gradient Ricci almost soliton \ref{eq1} with $h$ as a potential function, where the functions $h$ and $\lambda$ are given by
	
          \begin{equation} \label{eq6}
             \left\{ \begin{array}{lc}
               h(\xi)=\displaystyle\int \left[ c-(n-2)\int\psi\psi^{''}d\xi \right] \dfrac{1}{\psi^{2}}d\xi +k; \\ \\
                    
               \lambda(\xi)=\varepsilon_{i_{0}}\left\lbrace  \left[ \psi\psi^{''}-(n-1) (\psi^{'})^{2} \right] -  \dfrac{\psi^{'}}{\psi} \left[ c-(n-2) \displaystyle\int\psi\psi^{''}d\xi \right]\right\rbrace,    
            \end{array} \right. 
	     \end{equation}
	   	
   where $c$ and $k$ are constants.
\end{coro}

\begin{coro} \label{coro1.4}
    	If $(\mathbb{R}^{n}, g)$ is the Euclidean space, $F$ a Ricci-flat complete Riemannian manifold (if it is the case) and $0 < |\psi (x) |\leq c$ for some constant $c$, then the metrics in Theorem \ref{coro1.3}, Corollary \ref{coro1.1} and Corollary \ref{coro1.2} are complete.      	
\end{coro}

As a consequence of Corollary \ref{coro1.4} we obtain the following examples.

\begin{exam}\label{aaaaa}
 Considering $\alpha_{1}=\ldots=\alpha_{n-1}$, $\alpha_{n}=1$, $\psi(x_{1},\ldots,x_{n})=x_{n}$ and $\mathbb{R}^{n^{\ast}}_{+}=\{(x_{1},\ldots,x_{n})\in\mathbb{R}^{n}; x_{n}>0\}$, i.e., $\left(\mathbb{R}^{n^{\ast}}_{+},g_{can}=\frac{\delta_{ij}}{x_{n}^{2}}\right)=(\mathbb{H}^{n},g_{can})$, where $g_{can}$ is the standard metric of the hyperbolic space. Therefore, from Theorem \ref{coro1.3} we have that the product manifold $\mathbb{H}^{n}\times F^{m}$, in which $(F^{m}, g_{F})$ is a complete Ricci-flat manifold, is a complete gradient Ricci almost soliton with metric tensor $$ds^{2}= g_{can}+\frac{1}{x_{n}^{2}}g_{F},$$    
where the potential function is $$h(x_{1},\ldots,x_{m+n})=k-\frac{c}{x_{n}},$$ and $$\lambda(x_{1},\ldots,x_{m+n})=-\left[(m+n-1)+\frac{c}{x_{n}}\right].$$
\end{exam}

\begin{exam}
Consider $\alpha_{1}=\ldots=\alpha_{n-1}$, $\alpha_{n}=1$ and $\mathbb{R}^{n^{\ast}}_{+}=\{(x_{1},\ldots,x_{n})\in\mathbb{R}^{n}; x_{n}>0\}$, i.e., $\left(\mathbb{R}^{n^{\ast}}_{+},g_{can}=\frac{\delta_{ij}}{x_{n}^{2}}\right)=(\mathbb{H}^{n},g_{can})$, where $g_{can}$ is the standard metric of the hyperbolic space. Now, taking $\psi(x_{1},\ldots,x_{n})=x_{n}\nu$, where $\nu$ is any bounded and positive smooth function (cf. Corollary \ref{coro1.4})
we have that the product manifold $\mathbb{H}^{n}\times F^{m}$, in which $(F^{m}, g_{F})$ is a complete Ricci-flat manifold, is a complete gradient Ricci almost soliton with metric tensor $$ds^{2}= \frac{1}{\nu^{2}}\left(g_{can}+\frac{1}{x_{n}^{2}}g_{F}\right),$$    
where the functions $\lambda$ and $h$ are given by \ref{eq8}. Therefore, there exists many complete Ricci almost solitons conformal to Example \ref{aaaaa}. 
\end{exam}
	
\begin{exam}
If $\psi(x_{1},\ldots,x_{n})= e^{-x_{1}^{2}-\ldots-x_{n}^{2}}$ in Corollary \ref{coro1.1}, then $\mathbb{R}^{n}$ with Riemannian metric $$ds^{2}=e^{2(x_{1}^{2}+\ldots+x_{n}^{2})}(dx_{1}^{2}+\ldots+dx_{n}^{2})$$ is a complete gradient Ricci almost soliton   
    where the potential function is $$h(x_{1},\ldots,x_{n})=\frac{c}{2}e^{2x_{1}^{2}+\ldots+2x_{n}^{2}}+\dfrac{(n-2)}{2}(x_{1}^{2}+\ldots+x_{n}^{2})+k,$$ and 
    \begin{eqnarray*}
 \lambda(x_{1},\ldots,x_{n})&=&-[2(n-2)(x_{1}^{2}+\ldots+x_{n}^{2})+(3n-2)]e^{-2x_{1}^{2}-\ldots-2x_{n}^{2}}\\
 &+&2c[2(x_{1}^{2}+\ldots+x_{n}^{2})+1].
 \end{eqnarray*}
 \end{exam}

\begin{exam} 
Choosing n=2, $\alpha_{1}=\alpha_{2}=1$ and $\psi(x_{1},x_{2})=\frac{1}{1+(x_{1}+x_{2})^{2}}$, then from Corollary \ref{coro1.2} we have that $\mathbb{R}^{2}$ with metric $$ds^{2}=[1+2(x_{1}+x_{2})^{2}+(x_{1}+x_{2})^{4}](dx_{1}^{2}+dx_{2}^{2})$$ is a complete gradient Ricci almost soliton,   
	where the potential function is $$h(x_{1},x_{2})=c\left[(x_{1}+x_{2})+\frac{2}{3}(x_{1}+x_{2})^{3}+\frac{1}{5}(x_{1}+x_{2})^{5}\right]+k$$ and $$\lambda(x_{1},x_{2})=\frac{4(x_{1}+x_{2})^{2}}{[1+(x_{1}+x_{2})^{2}]^{4}}-\frac{2}{[1+(x_{1}+x_{2})^{2}]^{3}}+\frac{2c(x_{1}+x_{2})}{1+(x_{1}+x_{2})^{2}}.$$
\end{exam}


\section{Proofs of the Main Results}

\n \textbf{Proof of Theorem \ref{teo1.3}:} 

Let $(\mathbb{R}^{n}, g)$ be a pseudo-Euclidean space $n\geq3$ with coordinates $x=(x_{1}, ..., x_{n})$,  $(M, \tilde{g} )=(\mathbb{R}^{n},\bar{g})\times_{f} (F^{m}, g_{F})$ a warped product where $\tilde{g}= \overline{g} +f^{2}g_{F},\ \overline{g}=\frac{1}{\psi^{2^{}}}g, \ g_{ij}=\delta_{ij}\varepsilon_{i}$, $F$ is a pseudo-Riemannian Einstein manifold with constant Ricci curvature $\lambda_{F}$. Considering $X_{1},...,X_{n}\ \in \ \pounds (\mathbb{R}^{n})$  and $Y_{1},...,Y_{m}\ \in \ \pounds (F)$, where $\pounds (\mathbb{R}^{n})$ and $\pounds (F)$ are, respectively, the spaces of lifts of vector fields on $\mathbb{R}^{n}$ and $F$ to $\mathbb{R}^{n}\times_{f} F^{m}$.	
Since $\tilde{g}\left(Y_{i}, Y_{j} \right)=f^{2} g_{F}\left(Y_{i}, Y_{j} \right)$, from the warped structure (see \cite{B,O'neil}) we obtain
                     \begin{equation} \label{eq22}
                          \left\{ \begin{array}{lll}
                               Ric_{\tilde{g}}\left(X_{i}, X_{j} \right)&=&  Ric_{\overline{g}}\left(X_{i}, X_{j} \right) -\dfrac{m}{f} Hess_{\overline{g}} f\left(X_{i}, X_{j} \right), \ \forall \  i, j=1,...,n;\\ \\
                     
                              Ric_{\tilde{g}}\left(X_{i}, Y_{j} \right)&=& 0, \ \forall \  i=1,...,n;\ \ e \ \ j=1,...,m \\ \\
                     
                              Ric_{\tilde{g}}\left(Y_{i}, Y_{j} \right)&=&  Ric_{g_{F}}\left(Y_{i}, Y_{j} \right)\\ 
                              &-& \left(f\Delta_{\overline{g}}f +(m-1)\left|\nabla_{\overline{g}}f \right|^{2}\right) g_{F}\left(Y_{i}, Y_{j} \right), \ \forall \  i, j=1,...,m.
                          \end{array} \right.
                     \end{equation}	
It is well known (cf. \cite{B}), if $\overline{g}=\frac{1}{\psi^{2^{}}}g $, that
\begin{equation} \label{eq9}
                              Ric_{\overline{g}}=\dfrac{1}{\psi^{2}}\left\lbrace \left( n-2\right)\psi Hess_{g}\psi + \left[ \psi\Delta_{g}\psi-\left( n-1\right)\left|\nabla_{g}\psi \right|^{2}  \right]g\right\rbrace. 
                        \end{equation}
                        Furthermore, we get from the metric $g$
               $$\left( Hess_{g}\psi\right)\left(X_{i}, X_{j} \right)=\psi_{, x_{i}x_{j}}, \ \ \Delta_{g}\psi=\sum_{k}\varepsilon_{k}\psi_{, x_{k}x_{k}}, \ \ \left|\nabla_{\overline{g}}\psi \right|^{2}=\sum_{k}\varepsilon_{k} \left( \psi_{, x_{k}}\right) ^{2}$$
 
 inserting these expressions into \ref{eq9}, we get
               \begin{equation} \label{eq24}
                  \left\{ \begin{array}{lll}
                        Ric_{\overline{g}}\left(X_{i}, X_{j} \right)&=&   \left( n-2\right) \dfrac{\psi_{, x_{i}x_{j}}}{\psi}, \ \forall \  i\neq j=1,...,n;\\ \\
               
                       Ric_{\overline{g}}\left(X_{i}, X_{i} \right)&=& \dfrac{1}{\psi^{2}}\big\{ \left( n-2\right)\psi\psi_{, x_{i}x_{i}}  \\ 
                       &+& \big[ \psi\displaystyle\sum_{k}\varepsilon_{k}\psi_{, x_{k}x_{k}}-\left( n-1\right)\displaystyle\sum_{k}\varepsilon_{k} \psi_{, x_{k}}^{2}  \big]\varepsilon_{i} \big\}   , \ \forall \  i=1,...,n.
                 \end{array} \right.
               \end{equation}

Considering $X_{1},...,X_{n}\ \in \ \pounds (\mathbb{R}^{n})$ and \ref{eq1}, we have
                    $$ Ric_{\tilde{g}}\left(X_{i}, X_{j} \right)= \lambda\tilde{g}\left(X_{i}, X_{j} \right)-Hess_{\tilde{g}}\left( h\right) \left(X_{i}, X_{j} \right).$$
Then, from the first equation of \ref{eq22} 
                     \begin{eqnarray} \label{sia}Ric_{\overline{g}}\left(X_{i}, X_{j} \right)-\dfrac{m}{f}Hess_{\overline{g}} (f) \left(X_{i}, X_{j} \right)= \lambda\overline{g}\left(X_{i}, X_{j} \right)- Hess_{\overline{g}}(h) \left(X_{i}, X_{j} \right).
                     \end{eqnarray}

                     Let now $ \overline{g}=\frac{1}{\psi^{2^{}}}g $ be a conformal metric of $g$, $g_{ij}=\varepsilon_
{i}\delta_{ij}$, then the expressions of the Christoffel symbols are given by 

                 $$ \overline{\Gamma}_{ij}^{k}=0, \ \ \ \overline{\Gamma}_{ij}^{i}=-\dfrac{ \psi_{,x_{j}}}{\psi}, \ \ \ \overline{\Gamma}_{ii}^{k}= \varepsilon_{i} \varepsilon_{k}\dfrac{ \psi_{,x_{k}}}{\psi}\quad\mbox{and}\quad   \overline{\Gamma}_{ii}^{i}=-\dfrac{ \psi_{,x_{i}}}{\psi}.$$
                         
Then, for any smooth function $W$, the Hessian operator is given by                              
                  \begin{equation} \label{eq17}
                        \left\{ \begin{array}{lll}
                              Hess_{\overline{g}}(W)_{ij}=W_{, x_{i} x_{j}}+\dfrac{ \psi_{,x_{j}}}{\psi}W_{,x_{i}}+\dfrac{ \psi_{,x_{i}}}{\psi}W_{,x_{j}}, \ \ i\neq j;\\ \\
                         
                               Hess_{\overline{g}}(W)_{ii}=W_{, x_{i} x_{i}}+2\dfrac{ \psi_{,x_{i}}}{\psi}W_{,x_{i}}-\varepsilon_{i}\displaystyle\sum_{k=1}^{n}\varepsilon_{k} \dfrac{ \psi_{,x_{k}}}{\psi}W_{,x_{k}}.
                         \end{array} \right.
                   \end{equation}

Applying \ref{eq24} and the expression of Hessian in the metric $\bar{g}$ given by \ref{eq17} to \ref{sia}, we obtain
                     \begin{equation} \label{eq26}
                         \begin{array}{lll}
                             (n-2)f \psi_{,x_{i}x_{j}} &+& f \psi h_{,x_{i}x_{j}} -m \psi f_{,x_{i}x_{j}} - m\psi_{,x_{i}} f_{,x_{j}}\\
                           
                           &-&m\psi_{,x_{j}} f_{,x_{i}} + f \psi_{,x_{i}} h_{,x_{j}} + f \psi_{,x_{j}} h_{,x_{i}}=0,\quad 
                           1\leq i\neq j\leq n,
                        \end{array}
                     \end{equation}
                          and for all $i$           
                     \begin{equation} \label{eq27}
                         \begin{array}{lll}
                              \psi \left[ (n-2) f \psi_{,x_{i}x_{i}}+ f \psi h_{,x_{i}x_{i}} -m \psi f_{,x_{i}x_{i}} -2m \psi_{,x_{i}} f_{,x_{i}} +2f \psi_{,x_{i}} h_{,x_{i}} \right] + \\
                     
                              \varepsilon_{i} \sum^{n}_{k=1} \varepsilon_{k} \left[ f \psi \psi_{,x_{k}x_{k}} -(n-1) f \psi_{,x_{k}}^{2} +m \psi \psi_{,x_{k}} f_{,x_{k}} -f \psi \psi_{,x_{k}} h_{,x_{k}}\right]=\varepsilon_{i} \lambda f. \ \ \  
                        \end{array} 	
                     \end{equation}

Now, considering $Y_{1},...,Y_{m}\ \in \ \pounds (F)$, from \ref{eq1} and the third equation of \ref{eq22} we get 
                  \begin{equation} \label{eq28}
                  \begin{array}{lll}
                         Ric_{g_{F}}\left(Y_{i}, Y_{j} \right) &-& \left(f\Delta_{\overline{g}}f +(m-1)\left|\nabla_{\overline{g}}f \right|^{2}\right) g_{F}\left(Y_{i}, Y_{j} \right)\\
                         &-&\lambda f^{2}g_{F}\left(Y_{i}, Y_{j} \right) +\left( Hess_{\tilde{g}}h\right) \left(Y_{i}, Y_{j} \right)=0.
                         \end{array}
                  \end{equation}
It is a straightforward computation that
                    \begin{equation} \label{eq29}
                       \begin{array}{lll}
                           |\nabla_{\overline{g}}f|^{2}=\psi^{2}\Sigma_{k}\varepsilon_{k}f_{,x_{k}}^{2}\quad\mbox{and}\\ \\
                        
                           \triangle_{\overline{g}}f=\psi^{2}\Sigma_{k}\varepsilon_{k}f_{,x_{k}x_{k}}-(n-2)\psi\Sigma_{k}\varepsilon_{k}\psi_{, x_{k}}f_{,x_{k}}.
                       \end{array} 	
                  \end{equation}
 We also have that $F$ is an Einstein Manifold, thus            
                  \begin{equation} \label{eq30}
                         Ric_{g_{F}}\left(Y_{i}, Y_{j} \right)=\lambda_{F}g_{F}\left(Y_{i}, Y_{j} \right).
                  \end{equation}
Moreover,
\begin{eqnarray}\label{dan}
                   Hess_{\tilde{g}}(h)(Y_{i},Y_{j})&=&(Y_{i}Y_{j})(h)-(\nabla_{\tilde{g}Y_{i}}Y_{j})(h)\nonumber\\
                   &=&\left(\frac{\nabla_{\tilde{g}}f}{f}\right)(h)\tilde{g}(Y_{i},\,Y_{j})\nonumber\\
                   &=&f\bar{g}(\nabla_{\bar{g}}h,\nabla_{\bar{g}} f)g_{F}(Y_{i},Y_{j})\nonumber\\
                   &=&\left(f\psi^{2}\displaystyle\sum_{k}\varepsilon_{k}f_{,x_{k}}h_{,x_{k}}\right)g_{F}(Y_{i},Y_{j}).
                   \end{eqnarray}
Then, replacing this expression and \ref{eq29}, \ref{eq30}, \ref{dan} in \ref{eq28}, we get
                    \begin{equation} \label{eq31}
                          \sum^{n}_{k=1} \varepsilon_{k} \left[-f \psi^{2} f_{,x_{k}x_{k}} +(n-2) f \psi f_{,x_{k}} \psi_{,x_{k}} - (m-1) \psi^{2} f_{,x_{k}}^{2}+f \psi^{2} f_{,x_{k}}  h_{,x_{k}} \right]  = \lambda f^{2} 	-\lambda_{F}.
                    \end{equation}  
We assume that $f(\xi)$ and $\psi(\xi)$ are functions of $\xi$, where $\xi=\sum^{n}_{i=1}\alpha_{i} x_{i}$. Hence, we have
           $$\psi_{, x_{i}} = \alpha_{i} \psi^{'}, \ \ \ \ \ \ \psi_{, x_{i} x_{j}} = \alpha_{i} \alpha_{j} \psi^{''},  \ \ \ \ \ \    \psi_{, x_{i} x_{i}} = \alpha_{i}^{2} \psi^{''},$$ 	
           $$f_{, x_{i}} = \alpha_{i} f^{'}, \ \ \ \ \ \  f_{, x_{i} x_{j}} = \alpha_{i} \alpha_{j} f^{''}, \ \ \ \ \ \  f_{, x_{i} x_{i}} = \alpha_{i}^{2} f^{''} ,$$ 
and		
          $$|\nabla_{g}\psi|^{2} =  \left( \displaystyle\Sigma_{i=1}^{n} \varepsilon_{i} \alpha_{i}^{2} \right) (\psi^{'})^{2} = \varepsilon_{i_{0}} (\psi^{'})^{2}, \ \ \ \ \ \  \triangle_{g}\psi= \left(\displaystyle \Sigma_{i=1}^{n}\varepsilon_{i} \alpha_{i}^{2} \right) \psi^{''}=\varepsilon_{i_{0}} \psi^{''},$$ where $\displaystyle\sum_{i=1}^{n}\varepsilon_{i}\alpha_{i}^{2}=\varepsilon_{i_{0}}.$
          
We replace these expressions in \ref{eq26}, \ref{eq27} and \ref{eq31}.
From \ref{eq26}, we have
             $$\alpha_{i} \alpha_{j} \left\lbrace f \left[ (n-2) \psi^{''}+2\psi^{'} h^{'}+\psi h^{''}\right] -m\psi f^{''}-2m\psi^{'} f^{'}\right\rbrace =0, \ \ \ \forall i\neq j. $$
If there exists $i\neq j$ such that $\alpha_{i} \alpha_{j}\neq 0$, then this equation reduces to  the first equation of \ref{eq7}.
Likewise, from \ref{eq27}, we get
\begin{equation*}
\begin{array}{lll}
               &&\psi \alpha_{i}^{2}\left\lbrace  \left[ (n-2) \psi^{''}+2\psi^{'} h^{'}+\psi h^{''}\right] -m\psi f^{''}-2m\psi^{'} f^{'}\right\rbrace \\
               &+&  \varepsilon_{i}\Sigma_{k=1}^{n} \varepsilon_{k} \alpha_{k}^{2}\left[ f\psi\psi^{''}-(n-1)f \left( \psi^{'}\right)^{2}+m\psi\psi^{'} f^{'} - f \psi\psi^{'} h^{'} \right]=\varepsilon_{i}\lambda f,
               \end{array}
               \end{equation*}
applying the first equation of \ref{eq7} and taking $\sum_{k=1}^{n} \varepsilon_{k} \alpha_{k}^{2}=\varepsilon_{i_{0}}\neq 0$, the equation above reduces exactly to the second equation of \ref{eq7}. Finally, we have the third equation of \ref{eq7} by \ref{eq31}. The inverse statement of this theorem it is a straightforward computation.

 \vspace{.2in}

 \n \textbf{Proof of Theorem \ref{coro1.3}:} 
 
 We can rewrite the first equation of \ref{eq7} as
 
            $$f \psi h^{''}+2f \psi^{'} h^{'}+ (n-2)f\psi^{''}-m\psi f^{''}-2m \psi^{'} f^{'}=0.$$
 
 Making the change of variable $y=h^{'}$, the above equation becomes a first-order linear differential equation in $y$, 
 
          $$y^{'}=-2\dfrac{\psi^{'}}{\psi}y+\left[ -(n-2)\dfrac{\psi^{''}}{\psi}+m\dfrac{f^{''}}{f}+2m\dfrac{f^{'} \psi^{'}}{f \psi}\right].$$
 
Considering $\bar{f}(\xi)= -2\dfrac{\psi^{'}}{\psi}$ and $\bar{g}(\xi)=  -(n-2)\dfrac{\psi^{''}}{\psi}+m\dfrac{f^{''}}{f}+2m\dfrac{f^{'} \psi^{'}}{f \psi}$, by the linearity of the differential equations, we get
 
          $$y= \left[c+\int \bar{g}( \xi) e^{-\int \bar{f}(\xi) d \xi }  d \xi \right] e^{\int \bar{f}(\xi) d \xi}  d \xi.$$
 
Notice that
          $$\int \bar{f}(\xi) d \xi=-ln \psi^{2}.$$
Thus, 
         \begin{equation} \label{eq33}              
                 h(\xi)= k+ \int \left\lbrace c+ \int\left[ m\psi^{2}\dfrac{f^{''}}{f} -(n-2) \psi\psi^{''} +2m\dfrac{\psi\psi{'} f^{'}}{f} \right] d\xi \right\rbrace \frac{1}{\psi^{2}} d\xi.               
         \end{equation}
         
Hereafter, assume $f\psi=1$. Hence, 
                \begin{equation} \label{eq34}
     \psi f^{''} +2f^{'} \psi^{'}=-f \psi^{''}.                       
               \end{equation}
Note that we can rewrite equation \ref{eq33} using the equation \ref{eq34}
\begin{eqnarray} \label{eq35}      
                    h(\xi) &=& k+ \displaystyle \int \left\lbrace c+ \displaystyle \int\left[   m\psi^{2}\dfrac{f^{''}}{f} -(n-2) \psi\psi^{''} +2m\dfrac{\psi\psi^{'} f^{'}}{f} \right] d\xi \right\rbrace \frac{1}{\psi^{2}} d\xi  \nonumber\\ \nonumber\\
                           &=&  k+ \displaystyle \int \left\lbrace c+ \displaystyle \int \frac{\psi}{f} \left[ m\underbrace{\left( \psi f^{''}+2f^{'} \psi^{'}\right)}_{-f \psi^{''}} -(n-2)f \psi^{''} \right] d\xi \right\rbrace \frac{1}{\psi^{2}} d\xi \nonumber \\ \nonumber\\     
                    &=& k+ \displaystyle \int \left\lbrace c-(m+n-2) \displaystyle \int \psi \psi^{''}  d\xi \right\rbrace \frac{1}{\psi^{2}} d\xi.
               \end{eqnarray}

 Now, replacing $h$ and $h^{'}$ in the second equation of \ref{eq7} by the expressions given by \ref{eq35}, we obtain
\begin{equation} \label{eq37}              
             \lambda= \varepsilon_{i_{0}} \left\lbrace  \psi\psi^{''}-(n-1)(\psi^{'})^{2}+m\psi \psi^{'} \frac{f^{'}}{f} - c\dfrac{\psi^{'}}{\psi} +(m+n-2)\dfrac{\psi^{'}}{\psi} \displaystyle \int \psi \psi^{''} d\xi  \right\rbrace.                   
         \end{equation}
         
Applying the expressions $f=\dfrac{1}{\psi}$ and $f^{'}=-\dfrac{\psi^{'}}{\psi^{2}}$ to \ref{eq37}, we obtain 
\begin{eqnarray} \label{eq38}     
                  \lambda (\xi)=  \varepsilon_{i_{0}} \left\lbrace  \psi\psi^{''}-(m+n-1)(\psi^{'})^{2} -c \dfrac{\psi^{'}}{\psi} +(m+n-2)\dfrac{\psi^{'}}{\psi} \displaystyle \int \psi \psi^{''} d\xi  \right\rbrace .
              \end{eqnarray}
 
Now, we prove that the expression \ref{eq38} is well defined. 
The third equation of the system \ref{eq8} is given by
  	             	               
 	                $$ \varepsilon_{i_{0}} \left[-f \psi^{2} f^{''}+(n-2)f \psi f^{'}\psi^{'}-(m-1) \psi^{2} \left( f^{'}\right)^{2}+f \psi^{2}f^{'}h^{'} \right]  = \lambda f^{2}-\lambda_{F}. $$ 
Inserting $f=\dfrac{1}{\psi}, \ f^{'}=-\dfrac{\psi^{'}}{\psi^{2}},\ f^{''}=2\dfrac{\left( \psi^{'}\right) ^{2}}{\psi^{3}}-\dfrac{\psi^{''}}{\psi^{2}}$ and $\lambda_F =0$ into the above equation we obtain
\begin{eqnarray*}
\lambda (\xi) = \varepsilon_{i_{0}} \left[ \psi \psi^{''}-(m+n-1)(\psi^{'})^{2}-   \psi \psi^{'} h^{'} \right].
                     \end{eqnarray*}
                     Then, using the expression of $h^{'}$ given by equation \ref{eq35}, we get 
                     \begin{eqnarray}\label{eq39}
                     \lambda (\xi)=\varepsilon_{i_{0}} \left\lbrace  \psi\psi^{''}-(m+n-1)(\psi^{'})^{2} - c\dfrac{\psi^{'}}{\psi} +(m+n-2)\dfrac{\psi^{'}}{\psi} \displaystyle \int \psi \psi^{''} d\xi  \right\rbrace.
               \end{eqnarray}
 
Thus, equation \ref{eq39} is exactly the same as \ref{eq38}.
 
 \hfill $\Box$
 
 \vspace{.2in}

 \n \textbf{Proof of Theorem \ref{teo1.1}:} 

Let $(\mathbb{R}^{n}, g)$ be a pseudo-Euclidean space with coordinates $x=(x_{1}, ..., x_{n}),  \ g_{ij}=\varepsilon_{i}\delta_{ij}$. Since the metric $\overline{g}$ is a gradient Ricci almost soliton \ref{eq1}, we have
                    \begin{equation} \label{eq10}
                        Ric_{\overline{g}}+ Hess_{\overline{g}}(h)=\lambda\overline{g}, \ \ \ \lambda \in C^{\infty} (\mathbb{R}^{n}).
	\end{equation}

Now applying \ref{eq9} and \ref{eq17} to \ref{eq10}, we get
                  \begin{equation} \label{eq13}
                        \left\{ \begin{array}{llc}
                           (n-2) \psi_{,x_{i}x_{j}}+ \psi_{,x_{_{i}}}h_{,x_{j}}+{\psi_{,x_{_{j}}}h_{,x_{i}}+\psi h_{,x_{i}x_{j}}}=0 ,\quad i\neq j;\\ \\
                  
                            \psi[(n-2)\psi_{,x_{i}x_{i}}+ h_{,x_{i}x_{i}}\psi+2\psi_{,x_{i}}h_{,x_{i}}]\\
+ \varepsilon_{i}\big[ \psi \Delta_{g} \psi -(n-1)|\nabla_{g} \psi|^{2}-\psi\displaystyle\sum^{n}_{k=1}\varepsilon_{k}h_{,x_{k}}\psi_{,x_{k}}\big]= \varepsilon_{i}\lambda,\quad\mbox{for all}\quad i.
                       \end{array} \right.
                  \end{equation}	
Considering that $f(r)$ and $\psi(r)$ are functions of $r$, where  $r=\sum^{n}_{i=1}\varepsilon_{i}x_{i}^{2}$. Hence, we have	
$$\psi_{, x_{i}} = 2\varepsilon_{i}x_{i} \psi^{'}, \ \ \ \ \ \ \psi_{, x_{i} x_{j}} = 4\varepsilon_{i}\varepsilon_{j}x_{i} x_{j} \psi^{''},  \ \ \ \ \ \    \psi_{,x_{i} x_{i}} = 4x_{i}^{2} \psi^{''} +  2\varepsilon_{i}\psi^{'},$$			
           $$h_{,x_{i}} = 2\varepsilon_{i}x_{i}h^{'}, \ \ \ \ \ \  h_{, x_{i} x_{j}} = 4\varepsilon_{i}\varepsilon_{j}x_{i}x_{j}h^{''}, \ \ \ \ \ \  h_{,x_{i} x_{i}} = 4x_{i}^{2}h^{''} +  2\varepsilon_{i}h^{'}.$$ 

Furthermore,
           $$|\nabla_{g}\psi|^{2} = 4r(\psi^{'})^{2}, \ \ \ \ \ \  \triangle_{g}\psi=4r \psi^{''}+2n\psi^{'}.$$ 

Replacing these expressions in the first equation of \ref{eq13}, we get

         $$4\varepsilon_{i}\varepsilon_{j}x_{i} x_{j} \left[ (n-2) \psi^{''}+2\psi^{'} h^{'}+\psi h^{''}\right]=0, \ \ \ \forall i\neq j. $$

Therefore,
         \begin{equation} \label{eq14}	            	
             (n-2) \psi^{''}+2\psi^{'}h^{'}+\psi h^{''}=0 .		               
         \end{equation}

In an analogous way, from the second equation of \ref{eq13}, we get
\begin{eqnarray*}
        &&\varepsilon_{i}[4(n-1)\psi\psi^{'}+4r\psi\psi^{''}-4(n-1)r(\psi^{'})^{2}+2\psi^{2}h^{'}-4r\psi\psi^{'} h^{'}]\\
        &&+4x_{i}^{2}\psi \left[ (n-2) \psi^{''}+2\psi^{'} h^{'}+\psi h^{''} \right] =\varepsilon_{i}\lambda.
\end{eqnarray*}

Now applying \ref{eq14} to the above equation, it yields
\begin{equation} \label{eq15}		           
              4(n-1)\psi\psi^{'}+4r\psi\psi^{''}-4(n-1)r(\psi^{'})^{2}+2\psi^{2}h^{'}-4r\psi\psi^{'}h^{'}=\lambda.             
       \end{equation}

The reciprocal of this theorem can be easily verified.

\hfill $\Box$

\vspace{.2in}

\n \textbf{Proof of Corollary \ref{coro1.1}:} 

Note that we can rewrite \ref{eq14} as

              $$ h^{''} +2\dfrac{\psi^{'}}{\psi}h^{'}+(n-2)\dfrac{\psi^{''}}{\psi}=0. $$

Taking $y=h^{'}$ in this last equation, it is equivalent to first-order linear differential equation

            $$y^{'} +2\dfrac{\psi^{'}}{\psi}y+(n-2)\dfrac{\psi^{''}}{\psi}=0 $$

Solving this ordinary differential equation we get

            $$y= \left[ c-(n-2)\int \psi\psi^{''} dr \right] \dfrac{1}{\psi^{2}} $$

where $c>0$. Thus, we obtain the first equation of \ref{eq4}. Moreover, applying the first equation of \ref{eq4} to \ref{eq15}, we have the second equation of \ref{eq4}.

\hfill $\Box$

\vspace{.2in}

\n \textbf{Proof of Theorem \ref{teo1.2}:}

From \ref{eq17}, \ref{eq10} gives

                    \begin{equation} \label{eq18}
                            \left\{ \begin{array}{lll}
                                &&\psi_{,x_{i}x_{j}}=-\dfrac{1}{n-2} \left( {\psi_{,x_{_{i}}}h_{,x_{j}}+{\psi_{,x_{_{j}}}h_{,x_{i}}+\psi h_{,x_{i}x_{j}}}} \right) ,\quad i\neq j;\\ \\
                    
                                &&(n-2)\psi\psi_{,x_{i}x_{i}}+ \left[ \psi \Delta_{g} \psi-(n-1)|\nabla_{g} \psi|^{2} \right] \varepsilon_{i}\\
                                &+&2\psi\psi_{,x_{i}}h_{,x_{i}}+\psi^{2}h_{,x_{i}x_{i}}-\psi \varepsilon_{i}\displaystyle\sum^{n}_{k=1} \varepsilon_{k}h_{,x_{k}}\psi_{,x_{k}}= \lambda\varepsilon_{i}.
                            \end{array} \right.
                      \end{equation}

From now on, $h(\xi)$ and $\psi(\xi)$ are functions of $\xi$, where $\xi=\displaystyle\sum^{n}_{i=1}\alpha_{i} x_{i}$ and $\varepsilon_{i_{0}}=\displaystyle\sum^{n}_{i=1}\varepsilon_{i}\alpha_{i}^{2}$. Hence, we have $$\psi_{, x_{i}} = \alpha_{i} \psi^{'}, \ \ \ \ \ \ \psi_{, x_{i} x_{j}} = \alpha_{i} \alpha_{j} \psi^{''},  \ \ \ \ \ \    \psi_{, x_{i} x_{i}} = \alpha_{i}^{2} \psi^{''},$$ 	
              $$h_{, x_{i}} = \alpha_{i} h^{'}, \ \ \ \ \ \  h_{, x_{i} x_{j}} = \alpha_{i} \alpha_{j} h^{''}, \ \ \ \ \ \  h_{, x_{i} x_{i}} = \alpha_{i}^{2}h^{''} .$$ 
 Moreover,	
              $$|\nabla_{g}\psi|^{2} =  \left( \displaystyle\sum_{i=1}^{n} \varepsilon_{i} \alpha_{i}^{2} \right) (\psi^{'})^{2} = \varepsilon_{i_{0}} (\psi^{'})^{2}, \ \ \ \ \ \  \triangle_{g}\psi= \left(\displaystyle \sum_{i=1}^{n}\varepsilon_{i} \alpha_{i}^{2} \right) \psi^{''}=\varepsilon_{i_{0}} \psi^{''}.$$ 
Replacing these expressions in the first equation of \ref{eq18}, we get
             $$\alpha_{i} \alpha_{j} \left[ (n-2) \psi^{''}+2\psi^{'} h^{'}+\psi h^{''}\right]=0, \ \ \ \forall\, i\neq j. $$
If there exists $i\neq j$ such that $\alpha_{i} \alpha_{j}\neq 0$, then we obtain
\begin{equation} \label{eq19}	            	
       (n-2) \psi^{''}+2\psi^{'}h^{'}+\psi h^{''}=0,		               
            \end{equation}
which is exactly the first equation of \ref{eq5}.

Likewise, considering the second equation of \ref{eq18}, we get 
  \begin{eqnarray}  \label{x1}       
  \psi \alpha_{i}^{2}\left[ (n-2) \psi^{''}+2\psi^{'} h^{'}+\psi h^{''} \right] 
  +\varepsilon_{i_{0}} \left[ \psi\psi^{''}-(n-1) (\psi^{'})^{2} \right]\varepsilon_{i}  -\varepsilon_{i_{0}} \psi\psi^{'}h^{'} \varepsilon_{i} =\lambda\varepsilon_{i}.
  \end{eqnarray}
From the first equation of \ref{eq5} and \ref{x1}, we obtain the second equation in \ref{eq5}. 

We still have to consider the case in which $\alpha_{i_{0}}=1$ and $\alpha_{i}\neq0$ for all $i\neq i_{0}$. Then, the first equation of \ref{eq18} is trivially satisfied for all $i\neq j$.
Considering the second equation of \ref{eq18} for $i\neq i_{0}$, we get $$\varepsilon_{i_{0}} ( \psi\psi^{''}-(n-1) (\psi^{'})^{2} - \psi\psi^{'}h^{'}) =\lambda,$$ and hence the second equation of \ref{eq6} is satisfied. Considering $i = i_{0}$ in the second equation of \ref{eq18}, we get that the first equation \ref{eq5} is also verified.

\hfill $\Box$

\vspace{.2in}

\n \textbf{Proof of Corollary \ref{coro1.2}:}

From the first equation of \ref{eq5} we have
\begin{equation*}
                     h^{'}(\xi)= \dfrac{1}{\psi^{2}}\left[c-(n-2)\int \psi\psi^{''} d\xi \right].
                \end{equation*}
 Therefore, applying the above equation to the second equation of \ref{eq5}, we obtain our result. 

\hfill $\Box$

\vspace{.2in}

\vspace{.2in}

\n \textbf{Proof of Corollary \ref{coro1.4}:} 

Consider the Euclidean space $(R^{n}, g)$ $n\geq3$ and a metric $\overline{g}$ given by Corollaries \ref{coro1.1} and \ref{coro1.3}. If $0<|\psi (x)|\leq c$, then the metric $\overline{g}$ is complete, since there exists a constant $k >0$ such that for any vector $v \in R^{n}$, $|v|_{\overline{g}}\geq k |v|.$ We have that $M=(R^{n},\bar{g})\times_{f} F^{m},$ is complete if, and only if, $ (R^{n},\bar{g})  $ and $F^m$ are complete (see \cite{O'neil}). Then, the metrics obtained in Corollaries \ref{coro1.1} and \ref{coro1.3} are complete.

\hfill $\Box$


\vspace{.2in}



\begin{thebibliography} {99}


 
 \bibitem{BBR} {A. Barros, R. Batista, E. Ribeiro Jr.} - \textit{Rigidity of gradient Ricci almost solitons},
 Illinois J. Math. 56 (2012), 1267â€“-1279.
 
 
 
 \bibitem{BGR} {A. Barros, J. N. Gomes, E. Ribeiro Jr.} - \textit{A note on rigidity of the almost Ricci soliton},
 Arch. Math. (Basel) 100 (2013), 481â€“-490.




\bibitem{BBGG} {W. Batat, M. Brozos--VÃ¡zquez, E. GarcÃ­a--RÃ­o, S. Gavino-FernÃ¡ndez} -
 \textit{Ricci solitons on Lorentzian manifolds with large isometry groups}, Bull.
 London Math. Soc., 43 (2011), no. 6, 1219--1227.
 
 
\bibitem{B} {A. L. Besse} - \textit{Einstein Manifolds},
Springer-Verlag, Berlin, 1987.


\bibitem{BK} {V. Borges, K. Tenenblat} - \textit{Ricci almost solitons on semi-Riemannian warped products},
arXiv:1709.04604 (2017). 




 \bibitem{BCGG} {M. Brozos--VÃ¡zquez, G. Calvaruso, E. GarcÃ­a--RÃ­o, S. Gavino-FernÃ¡ndez} -
\textit{Three dimension Lorentzian homogeneous Ricci solitons}, Israel. J. Math., 188 (2012), 385--403.




 \bibitem{BGG} {M. Brozos--VÃ¡zquez, E. GarcÃ­a--RÃ­o, S. Gavino-FernÃ¡ndez} -
\textit{Locally conformally flat Lorentzian gradient Ricci solitons}, J. Geom. Anal., 23 (2013), no. 3, 1196--1212.



\bibitem{BGRR} {M. Brozos-VÃ¡zquez, E. GarcÃ­a--RÃ­o, X. Valle-Regueiro} - \textit{Half conformally flat gradient Ricci almost solitons}, Proc. R. Soc. A 472.2189 (2016): 20160043.
 
 
  \bibitem{CEGGV}{ E. Calvi\~no-Louzao, M. Fern\'andez-L\'opez, E. Garc\'ia-R\'io, R. V\'azquez-Lorenzo} - \textit{ Homogeneous Ricci almost solitons}, Israel J. Math., 220 (2017), no. 2, 531--546.  
 
 
 
 \bibitem{CA} {G. Catino} - \textit{Generalized quasi-Einstein manifolds with harmonic Weyl tensor}, Math. Z., 271, no. 3--4, 751--756, 2012.


\bibitem{catino} {G. Catino, et al.} - \textit{The Ricci-Bourguignon flow}, Pacific. J. Math., 287.2 (2017): 337-370.



\bibitem{FFGP} {F. E. S. Feitosa, A. A. Freitas Filho, J. N. V. Gomes, R. S. Pina} - \textit{Gradient Ricci almost soliton warped product}, arXiv:1507.03038v2, (2017). 




\bibitem{O'neil}
{B. O'Neill} - \textit{Semi--Riemannian geometry with applications to relativity}, (Academic Press, New York), 1983.



\bibitem{Onda} {K. Onda} - \textit{Lorentzian Ricci solitons on $3$-dimensional Lie groups}, Geom. Dedicata,
 {\bf 147} (2010), 313--322.
 
 
\bibitem{PRRS} {S. Pigola, M. Rigoli, M. Rimoldi, A. Setti} - \textit{Ricci almost solitons}, Ann. Scuola Norm. Sup. Pisa Cl. Sci., (5) Vol. X (2011),757--799. 
 

 

 \bibitem{PK} {R. Pina, K. Tenenblat} - \textit{On solitons of the Ricci curvature equation and the Einstein equation}, Israel J. Math., 171 (2009), 61--76. 
 




 


\end{thebibliography}
\end{document}